\documentclass[12pt]{article}
\usepackage{amsmath}
\usepackage{amssymb}
\usepackage{amsthm}
\usepackage{amsfonts}
\usepackage{amscd}
\usepackage[mathscr]{eucal}
\newcommand{\Z} {{\mathbb  Z}}
\newcommand{\Q}{{\mathbb  Q}}

\newcommand{\R} {{\mathbb R}}

\begin{document}
\parindent  25pt
\baselineskip  10mm
\textwidth  15cm    \textheight  23cm \evensidemargin -0.06cm
\oddsidemargin -0.01cm

\title{{A graphical method to calculate Selmer groups of several
families of non-CM elliptic curves}}
\author{\mbox{}  {Fei Li and  Derong Qiu }
\thanks{ \quad E-mail:
derong@mail.cnu.edu.cn } \\
(School of Mathematical Sciences, \\
Institute of Mathematics and Interdisciplinary Science, \\
Capital Normal University, Beijing 100048, P.R.China ) }

\date{}
\maketitle
\parindent  24pt
\baselineskip  10mm
\parskip  0pt

\par     \vskip  0.2 cm

{\bf Abstract} \quad In this paper, we extend the ideas of Feng
[F1], Feng-Xiong [FX] and Faulkner-James [FJ] to calculate the
Selmer groups of elliptic curves $ y^{2} = x ( x + \varepsilon p D )
( x + \varepsilon q D ). $
\par   \vskip 0.2cm

\par  \vskip  0.1 cm


\par   \vskip 0.3 cm

\hspace{-0.6cm}{\bf 1 \ \ Introduction and Main Results}

\par \vskip 0.2 cm

In this paper, we consider the following elliptic curves \\
$ \begin{array}{l} E = E_{\varepsilon }: y^{2} = x ( x + \varepsilon
p D ) ( x + \varepsilon q D ),  \qquad
\qquad \quad (1.1) \\
E^{\prime } = E_{\varepsilon }^{\prime }: y^{2} = x^{3} - 2
\varepsilon ( p + q ) D x^{2} + 4^{m} D^{2} x, \qquad \qquad \quad
(1.2)
\end{array} $ \\
where $ \varepsilon = \pm 1, \ p $ and $ q $ are odd prime numbers
with $ q - p = 2^{m}, m \geq 1 $ and $ D = D_{1} \cdots D_{n} $ is a
square-free integer with distinct primes $ D_{1}, \cdots , D_{n}. $
Moreover, $ 2 \nmid D, \ p \nmid D $ and $ q \nmid D. $ For each $
D_{i}, $ denote $ \widehat{D_{i}} = D / D_{i} \ ( \widehat{D_{1}}= 1
\ \text{if } \ D = D_{1} ). $ We write $ E = E_{+}, E^{\prime} =
E^{\prime}_{+} $ if $ \varepsilon = 1, $ and $ E = E_{-}, E^{\prime}
= E^{\prime}_{-} $ if $ \varepsilon = -1. $
\par \vskip 0.1 cm
There is an isogeny $ \varphi $ of degree $ 2 $ between $ E $ and $
E^{\prime} $ as follows:
$$ \varphi : E \longrightarrow E^{\prime}, \quad ( x, y ) \longmapsto
( y^{2} / x^{2}, \ y( pqD^{2} - x^{2} ) / x^{2}). $$ The kernel is $
E [ \varphi ] = \{ O, \ (0, 0) \}, $ and the dual isogeny of $
\varphi $ is $$ \widehat{ \varphi } : E^{\prime } \longrightarrow E,
\quad ( x, y ) \longmapsto ( y^{2} /4 x^{2}, \ y( 4^{m} D^{2} -
x^{2} ) /8 x^{2}) $$ with kernel $ E^{\prime} [ \widehat{ \varphi }
] = \{ O, \ (0, 0) \} $ (see [S, p.74]).
\par \vskip 0.1 cm

In this paper, we extend the ideas of Feng [F1], Feng-Xiong [FX] and
Faulkner-James [FJ] to calculate the $ \varphi ( \widehat{\varphi
})-$Selmer groups $ S^{(\varphi )} (E / \Q ) $ and $
S^{(\widehat{\varphi })} (E^{\prime} / \Q ). $ The main results are
as follows:
\par \vskip 0.2 cm

{\bf Theorem 1.1} \ Let $ D = D_{1}D_{2} \cdots D_{s}D_{s+1} \cdots
D_{n} $ with $ \left(\frac{p q}{D_{i}} \right) = 1 \ (i \leq s) \ $
and $ \left(\frac{p q}{D_{j}}\right) = -1 \ ( s < j \leq n )$ for
some non-negative integer $ s \leq n. $ If $ m, p $ and $ D $
satisfy one of the following conditions: \\
(1) \ $ m = 1, \ p D \equiv 5 ( \text{mod} 8) $ and $ D \equiv 3 (
\text{mod} 4); $ \\
(2) \ $ m = 1, \ p D \equiv 1 ( \text{mod} 8) $ and $ D \equiv 1 (
\text{mod} 4); $ \ (3) \ $ m = 2; $ \\ (4) \ $ p D \equiv 3 (
\text{mod} 8); $ \ (5) \ $ m = 3, \ p D \equiv 1 ( \text{mod} 4); $
\ (6) \ $ m = 4, \ p D \equiv 7 ( \text{mod} 8), $ \\
then $ \sharp S^{(\varphi )} (E / \Q ) = \sharp \{(V_{1},
V_{2})\mapsto_{e}G(+D): \ -1, p, q, D_{k} \notin V_{1}; \ s < k \leq
n \}. $ In the other cases, $ \sharp S^{(\varphi )} (E / \Q ) =
\sharp \{(V_{1}, V_{2})\mapsto_{e}G(+D): \ -1 ,p , q , D_{k} \in
V_{2}; s < k \leq n \} + \sharp \{(V_{1}, V_{2}) \mapsto_{qe}G(+D):
\ -1, p, q, D_{k} \notin V_{1}; \ s < k \leq n \}. $ Here $ G(+D) $
is the directed graph (see the following Definition 2.5).
\par \vskip 0.2 cm

{\bf Theorem 1.2.} \ Let $ D = D_{1}D_{2} \cdots D_{s}D_{s+1} \cdots
D_{n} $ with $ \left( \frac{p q}{D_{i}} \right) = 1 \ (i \leq s) $
and $ \left(\frac{p q}{D_{j}}\right) = -1 \ (s < j \leq n) $ for
some non-negative integer $ s \leq n, $ then $ \sharp
S^{(\widehat{\varphi })} (E^{\prime } / \Q ) = 2 \sharp \{(V_{1},
V_{2})\mapsto_{e} g(+D): \ \pm 2 \not \in V_{1} \}. $ Here $ g(+D) $
is the directed graph (see the following Definition 2.8).
\par \vskip 0.2 cm

{\bf Theorem 1.3.} \ Let $ D = D_{1}D_{2} \cdots D_{s}D_{s+1} \cdots
D_{n} $ with $ \left(\frac{p q}{D_{i}} \right) = 1 \ (i \leq s) $
and $ \left(\frac{p q}{D_{j}}\right) = -1 \ (s < j \leq n) $ for
some non-negative integer $ s \leq n. $ If $ m, p $ and $ D $
satisfy one of the following conditions: \\
(1) \ $ m = 1, \ p D \equiv 5, 7 ( \text{mod} 8) $ and $ D \equiv 3
(\text{mod} 4); $ \\
(2) \ $ m = 1, \ p D \equiv \pm 3 ( \text{mod} 8) $ and $ D \equiv 1
(\text{mod} 4); $ \\ (3) \ $ m = 2; $  \ (4) \ $ m = 3, \ p D \not
\equiv 1 ( \text{mod} 8); $ \\ (5) \ $ m = 4, \ p D \equiv 1 (
\text{mod} 4), $ \ (6) \ $ m \geq 5,
\ p D \equiv 5 ( \text{mod} 8), $ \\
then $ \sharp S^{(\varphi )} (E / \Q ) = \sharp \{(V_{1},
V_{2})\mapsto_{e} G(-D): \ p, q, D_{k} \in V_{2}; s < k \leq n \}; $
\ In other cases, $ \sharp S^{(\varphi )} (E / \Q ) = \sharp
\{(V_{1}, V_{2}) \mapsto_{e} G(-D): \ p, q, D_{k} \in V_{2}; s < k
\leq n \} + \sharp \{(V_{1}, V_{2})\mapsto_{qe} G(-D): \ p, q, D_{k}
\in V_{2}; s < k \leq n \}. $ Here $ G(-D) $ is the directed graph
(see the following Definition 2.10).
\par  \vskip  0.2 cm

{\bf Theorem 1.4.} \ Let $ D = D_{1}D_{2} \cdots D_{s}D_{s+1} \cdots
D_{n} $ with $ \left(\frac{p q}{D_{i}}\right) = 1 \ (i \leq s) $ and
$ \left(\frac{p q}{D_{j}}\right) = -1 \ (s < j \leq n ) $ for some
non-negative integer $ s \leq n, $ then \ $ \sharp
S^{(\widehat{\varphi })} (E^{\prime } / \Q ) = 2 \sharp \{(V_{1},
V_{2})\mapsto_{e}g(-D): \ -1, \pm 2 \notin V_{1} \}. $ Here $ g(-D)
$ is the directed graph (see the following Definition 2.13).
\par \vskip 0.1 cm
Moreover, another result about the Selmer group of elliptic curves
in (1.1) for all integers $ m \geq 2 $ is given in the appendix.

\par \vskip 0.2 cm
\hspace{-0.6cm}{\bf 2 \ \ Proofs of Theorems}

\par \vskip 0.2 cm

 Let $ M_{\Q} $ be the set of all places of $ \Q, $ including
the infinite $ \infty . $ For each place $ p, $ denote by $ \Q_{p} $
the completion of $ \Q $ at $ p $, and if $ p $ is finite, denote by
$ v_{p} $ the corresponding normalized additive valuation, so $
v_{p}(p) = 1. $ Let $ S = \{ \infty , 2, p, q, D_{1}, \cdots , D_{n}
\}, $ and define a subgroup of $ \Q^{\star } / \Q^{\star ^{2}} $ as
follows: \ $ \Q(S, 2)= <-1> \times <2> \times <p> \times <q> \times
<D_{1}> \times \cdots \times <D_{n}> \cong \left(\Z / 2 \Z
\right)^{n + 4}. $ \ For any subset $ A \subset \Q^{\star }, $ we
write $ <A> $ for the subgroup of $ \Q^{\star } / \Q^{\star ^{2}} $
generated by all the elements in $ A. $ For each $ d \in \Q(S, 2), $
define the curves \\
$ C_{d} : \ d w^{2} = d^{2} - 2 \varepsilon ( p + q) D d z^{2} +
4^{m}D^{2} z^{4}, $ \ and \\
$ C_{d}^{\prime} : \ d w^{2} = d^{2} +  \varepsilon ( p + q) D d
z^{2} + p q D^{2} z^{4}. $
\par  \vskip 0.1 cm
We have the following propositions 2.1 $ \sim $ 2.4 in determining
the local solutions of these curves $ C_{d} $ and $ C_{d}^{\prime}.
$ The proofs are similar to those in [LQ], so we omit the details.
\par  \vskip 0.2 cm

{\bf Proposition 2.1 } \ We assume $ \varepsilon = 1 $ and the
elliptic curve $ E = E_{+} $ be as in (1.1).
\par  \vskip 0.1 cm
(A) \ For $ d \in \Q (S, 2), $ if one of the following
conditions holds: \\
(1) \ $ d < 0; $ \quad  (2) \ $ p \mid d; $ \quad  (3) \ $ q \mid
d. $ \\
Then $ d \notin S^{(\varphi )} (E / \Q ). $ \ Moreover, if $ d >
0, $ then $ C_{d}( \R ) \neq \emptyset . $
\par  \vskip 0.1 cm
(B) \ For each $ d > 0, 2\mid d \mid 2 D, d  \in \Q (S, 2), $ we
have
\par  \vskip 0.1 cm
(1) \ if \ $ m = 1, $ then $ C_{d} ( \Q _{2}) \neq \emptyset \
\Longleftrightarrow \ \frac{d}{2} - 2 D ( p + 1) + \frac{2D^{2}}{d}
\equiv 2 (\text{mod} 16 ); $ \\
if \ $ m = 2, $ then $ C_{d} ( \Q _{2}) = \emptyset ; $ \\
if \ $ m = 3, $ then $ C_{d} ( \Q _{2}) \neq \emptyset \
\Longleftrightarrow \ d - D ( p + 4) + \frac{4D^{2}}{d}
\equiv 1 (\text{mod} 8); $ \\
if \ $ m = 4, $ then $ C_{d} ( \Q _{2}) \neq \emptyset \
\Longleftrightarrow \ d - D p  \equiv 1 (\text{mod} 8 ); $ \\
if \ $ m \geq 5, $ then $ C_{d} ( \Q _{2}) \neq \emptyset \
\Longleftrightarrow \ D p \equiv 7 (\text{mod} 8) $ \ or \ $ d - D p
\equiv 1 (\text{mod} 8). $
\par  \vskip 0.1 cm
(2) \ For each odd prime number $ l \mid \frac{2 p q D}{d}, \ C_{d}
( \Q _{l}) \neq \emptyset \ \Longleftrightarrow \ \left(\frac{d}{l}
\right) = 1. $
\par  \vskip 0.1 cm
(3) \ For each odd prime number $ l \mid d , \ C_{d} ( \Q _{l}) \neq
\emptyset \ \Longleftrightarrow \ \left(\frac{p D d
l^{-2}}{l}\right) = \left(\frac{q D d l^{-2}}{l}\right) = 1. $
\par  \vskip 0.1 cm
(C) \ For $  d > 0,d \mid  D, d  \in \Q (S, 2), $ we have
\par  \vskip 0.1 cm
(1) \ if \ $ m = 1, $ then $ C_{d} ( \Q _{2}) \neq \emptyset \
\Longleftrightarrow \ d \equiv 1 (\text{mod} 4 ); $ \\
if \ $ m = 2, $ then $ C_{d} ( \Q _{2}) \neq \emptyset \
\Longleftrightarrow d \equiv 1 (\text{mod} 4 ) $ or $ 2 d -  D ( p +
2)  \equiv 1 (\text{mod} 8); $ \\
if \ $ m \geq 3, $ then $ C_{d} ( \Q _{2}) \neq \emptyset \
\Longleftrightarrow \ d \equiv 1 (\text{mod} 4) $ or $ d -  D  p
\equiv 0(\text{mod} 4 ). $
\par  \vskip 0.1 cm
(2) \ For each prime number $ l \mid \frac{p q D}{d}, \ C_{d} ( \Q
_{l}) \neq \emptyset \ \Longleftrightarrow \
\left(\frac{d}{l}\right) = 1. $
\par  \vskip 0.1 cm
(3) \ For each prime number $ l \mid d , \ C_{d} ( \Q _{l}) \neq
\emptyset \ \Longleftrightarrow \ \left(\frac{pdDl^{-2}}{l} \right)
= \left(\frac{qdDl^{-2}}{l}\right) = 1. $
\par  \vskip 0.2 cm

{\bf Proposition 2.2 } \ We assume $ \varepsilon = 1 $ and the
elliptic curve $ E^{\prime} = E^{\prime}_{+} $ be as in (1.2).
\par  \vskip 0.15 cm
(A) \ (1) \ For any $ d \in \Q (S, 2), \ C_{d}^{\prime } ( \R ) \neq
\emptyset . $ If $ 2 \mid d, $ then $ d \notin S^{(\widehat{\varphi
})} (E^{\prime } / \Q ). $ \\
(2) \ $ \{ 1, pq, -pD, -qD \} \subset S^{(\widehat{\varphi })}
(E^{\prime } / \Q ). $
\par  \vskip 0.1 cm

(B) \ For each $  \ d \in \Q (S, 2) $ satisfying $ d \mid p D, $ we
have
\par  \vskip 0.1 cm
(B1) \ (1) \ If $ m = 1, $ then $ C_{d}^{\prime} ( \Q _{2}) \neq
\emptyset $ \ if and only if one of the following conditions holds:
(a) \ $ \ d \equiv 1 (\text{mod} 8), $ \ (b) \ $ ( d + p D ) ( d + q
D ) \equiv 0 (\text{mod} 16), $ \ (c) \ $ \frac{p q
D^{2}}{d} \equiv 1 ( \text{mod} 8); $ \\
(2) \ If $ m = 2, $ then $ C_{d}^{\prime} ( \Q _{2}) \neq \emptyset
$ \ if and only if one of the following conditions holds: \\ (a) \ $
\ d \equiv 1 (\text{mod} 8), $ \ (b) \ $ \frac{p q D^{2}}{d} \equiv
1 ( \text{mod} 8), $ \\ (c) \ $  d + p D \equiv 0 ( \text{mod} 4), $
\ (d) \ $  d \equiv 3 ( \text{mod} 4) $ \ and \ $
( p + 2 )D \equiv 1 ( \text{mod} 8 ); $ \\
(3) \ If $ m = 3, $ then $ C_{d}^{\prime} ( \Q _{2}) \neq \emptyset
$ \ if and only if one of the following conditions holds: \\ (a) \ $
d \equiv 1 (\text{mod} 8), $ \ (b) \ $ \frac{p q D^{2}}{d} \equiv 1
( \text{mod} 8), $ \ (c) \ $  d + p D \equiv 0 ( \text{mod} 8), $ \\
(d) \ $  d \equiv 3 ( \text{mod} 4), $ \ and \ $ d + p D \equiv 4 (
\text{mod} 8), $ \ (e) \ $  d \equiv 5 ( \text{mod} 8) $ \ and
\ $ d + p D \equiv 2 ( \text{mod} 4); $  \\
(4) \ If $ m = 4, $ then $ C_{d}^{\prime} ( \Q _{2}) \neq \emptyset
$ \ if and only if one of the following conditions holds: \\ (a) \ $
d \equiv 1 (\text{mod} 8), $ \ (b) \  $ \frac{p q D^{2}}{d} \equiv 1
( \text{mod} 8), $ \ (c) \ $  d + p D \equiv 0 ( \text{mod} 8), $ \\
(d) \ $  d \equiv 1 ( \text{mod} 8) $ \ and \ $ d + p D \equiv 2 (
\text{mod} 4), $ \ (e) \ $  d \equiv 5 ( \text{mod} 8) $ \ and
\ $ d + p D \equiv 4 ( \text{mod} 8); $  \\
(5) \ If $ m \geq 5, $ then $ C_{d}^{\prime} ( \Q _{2}) \neq
\emptyset $ \ if and only if one of the following
conditions holds: \\
(a) \ $ d \equiv 1 (\text{mod} 8), $ \ (b) \ $ \frac{p q D^{2}}{d}
\equiv 1 ( \text{mod} 8), $ \ (c) \ $  d + p D \equiv 0 ( \text{mod}
8). $
\par  \vskip 0.1 cm
(B2) \ $ C_{d}^{\prime} ( \Q _{p}) \neq \emptyset $ and $
C_{d}^{\prime} ( \Q _{q}) \neq \emptyset . $
\par  \vskip 0.1 cm
(B3) \ For each prime $ l\mid D , l \nmid d ,\quad C_{d}^{\prime} (
\Q _{l}) \neq \emptyset \Longleftrightarrow \ \left( 1 -
\left(\frac{d}{l}\right) \right) \left( 1 - \left(\frac{p q
d}{l}\right) \right) = 0. $
\par  \vskip 0.1 cm
(B4) \ For each prime $ l \mid D , l \mid d, \ C_{d}^{\prime} ( \Q
_{l}) \neq \emptyset \Longleftrightarrow  \left( 1 - \left(
\frac{-pdDl^{-2}}{l}\right) \right) \left( 1 -
\left(\frac{-qdDl^{-2}}{l}\right) \right) = 0. $
\par  \vskip 0.2 cm

{\bf Proposition 2.3 } \ We assume $ \varepsilon = - 1 $ and the
elliptic curve $ E = E_{-} $ be as in (1.1).
\par  \vskip 0.1 cm
(A) \ For $ d \in \Q (S, 2), $ if one of the following conditions
holds: \\
(1) \ $ p \mid d; $ \quad  (2) \ $ q \mid d. $ \\
Then $ d \notin S^{(\varphi )} (E / \Q ). $ \ Moreover,  $ C_{d}( \R
) \neq \emptyset . $
\par  \vskip 0.1 cm
(B) \ For each $ 2\mid d ,d \mid 2 D, d  \in \Q (S, 2), $ we have \\
(1) \ if \ $ m = 1, $ then $ C_{d} ( \Q _{2}) \neq \emptyset \
\Longleftrightarrow \frac{d}{2} + 2 D ( p + 1) + \frac{2D^{2}}{d}
\equiv 2 (\text{mod} 16); $ \\
if \ $ m = 2, $ then $ C_{d} ( \Q _{2}) = \emptyset; \ $ \\
if \ $ m = 3, $ then $  C_{d} ( \Q _{2}) \neq \emptyset \
\Longleftrightarrow
\ d +  D ( p + 4) + \frac{4D^{2}}{d} \equiv 1 (\text{mod} 8 ); $ \\
if \ $ m = 4, $ then $ C_{d} ( \Q _{2}) \neq \emptyset \
\Longleftrightarrow \
d + D p  \equiv 1 (\text{mod} 8 ) ;$ \\
if \ $ m \geq 5, $ then $ C_{d} ( \Q _{2}) \neq \emptyset \
\Longleftrightarrow \ D p \equiv 1 (\text{mod} 8) $ \ or \ $ d + D p
\equiv 1 (\text{mod} 8). $
\par  \vskip 0.1 cm
(2) \ For each odd prime number $ l \mid \frac{2 p q D}{d}, \ C_{d}
( \Q _{l}) \neq \emptyset \ \Longleftrightarrow \
\left(\frac{d}{l}\right) = 1. $ \\
(3) \ For each odd prime number $ l \mid d , \ C_{d} ( \Q _{l}) \neq
\emptyset \ \Longleftrightarrow \ \left(\frac{-pdDl^{-2}}{l}\right)
= \left(\frac{-qdDl^{-2}}{l}\right) = 1. $
\par  \vskip 0.1 cm
(C) \ For $ d \mid  D, d  \in \Q (S, 2), $ we have
\par  \vskip 0.2 cm
(1) \ if \ $ m = 1, $ then $ C_{d} ( \Q _{2}) \neq \emptyset \
\Longleftrightarrow
\ d \equiv 1 (\text{mod} 4); $ \\
if \ $ m = 2, $ then $ C_{d} ( \Q _{2}) \neq \emptyset \
\Longleftrightarrow \ d \equiv 1 (\text{mod} 4) $ or $ 2 d +  D ( p
+ 2)  \equiv 1
(\text{mod} 8); $ \\
if \ $ m \geq 3, $ then $ C_{d} ( \Q _{2}) \neq \emptyset \
\Longleftrightarrow\ d \equiv 1 (\text{mod} 4 ) $ or $ d + D p
\equiv 0(\text{mod} 4). $
\par  \vskip 0.1 cm
(2) \ For each prime number $ l \mid \frac{p q D}{d}, \ C_{d} ( \Q
_{l}) \neq \emptyset \ \Longleftrightarrow \
\left(\frac{d}{l}\right) = 1. $
\par  \vskip 0.1 cm
(3) \ For each prime number $ l \mid d, \ C_{d} ( \Q _{l}) \neq
\emptyset \ \Longleftrightarrow \ \left(\frac{-pdDl^{-2}}{l}\right)
= \left(\frac{-qdDl^{-2}}{l}\right) = 1. $
\par  \vskip 0.2 cm

{\bf Proposition 2.4.} \ We assume $ \varepsilon = - 1 $ and the
elliptic curve $ E^{\prime} = E^{\prime}_{-} $ be as in (1.2).
\par  \vskip 0.1 cm
(A) \ (1) \ For any $ d \in \Q (S, 2) $ and $ d > 0, \ C_{d}^{\prime
} ( \R ) \neq \emptyset . $ If $ 2 \mid d $ or $ d < 0, $ then $ d
\notin
S^{(\widehat{\varphi })} (E^{\prime } / \Q ). $ \\
(2) \ $ \{ 1, pq, pD, qD \} \subset S^{(\widehat{\varphi })}
(E^{\prime } / \Q ). $
\par  \vskip 0.1 cm
(B) \ For each $  \ d \in \Q (S, 2), d \mid p D , d > 0 , $ we have \\
(B1) \ (1) If $ m = 1, $ then $ C_{d}^{\prime} ( \Q _{2}) \neq
\emptyset $ if and only if one of the following conditions holds: \\
(a) \ $ \ d \equiv 1 (\text{mod} 8), $ \ (b) \ $ ( d - p D ) ( d - q
D ) \equiv 0 (\text{mod} 16), $ \ (c) \ $
\frac{p q D^{2}}{d} \equiv 1 ( \text{mod} 8); $ \\
(2) \ If $ m =2, $ then $ C_{d}^{\prime} ( \Q _{2}) \neq \emptyset $
if and only if one of the following conditions holds: \\
(a) \ $ d \equiv 1 (\text{mod} 8), $  \ (b) \ $ \frac{p q D^{2}}{d}
\equiv 1 ( \text{mod} 8), $ \\
(c) \ $  d - p D \equiv 0 ( \text{mod} 4), $  \ (d) \ $ d \equiv 3 (
\text{mod} 4) $ \ and \ $ ( p + 2 ) D
 \equiv 7 ( \text{mod} 8); $ \\
(3) \ If $ m = 3, $ then $ C_{d}^{\prime} ( \Q _{2}) \neq \emptyset
$ \ if and only if one of the following conditions holds: \\
(a) \ $ d \equiv 1 (\text{mod} 8), $  \ (b) \ $ \frac{p q D^{2}}{d}
\equiv 1 ( \text{mod} 8), $ \ (c) \ $  d - p D \equiv 0 ( \text{mod}
8), $ \\
(d) \ $  d \equiv 3 ( \text{mod} 4) $ \ and \ $ d - p D \equiv 4 (
\text{mod} 8), $ \ (e) \ $ d \equiv 5 ( \text{mod} 8) $ \ and
\ $ d - p D \equiv 2 ( \text{mod} 4). $ \\
(4) \ If $ m = 4, $ then $ C_{d}^{\prime} ( \Q _{2}) \neq \emptyset
$ \ if and only if one of the following conditions holds: \\
(a) \ $ d \equiv 1 (\text{mod} 8), $ \ (b) \ $ \frac{p q D^{2}}{d}
\equiv 1 ( \text{mod} 8), $ \ (c) \ $  d - p D \equiv 0 ( \text{mod}
8), $ \\
(d) \ $ d \equiv 1 ( \text{mod} 8) $ \ and \ $ d - p D \equiv 2 (
\text{mod} 4), $ \ (e) \ $ d \equiv 5 ( \text{mod} 8) $ \ and
\ $ d - p D \equiv 4 ( \text{mod} 8); $ \\
(5) \ If $ m \geq 5, $ then $ C_{d}^{\prime} ( \Q _{2}) \neq
\emptyset $ \ if and only if one of the following conditions holds: \\
(a) \ $ d \equiv 1 (\text{mod} 8), $ \ (b) \ $ \frac{p q D^{2}}{d}
\equiv 1 ( \text{mod} 8), $ \ (c) \ $ d - p D \equiv 0 ( \text{mod}
8). $
\par  \vskip 0.1 cm
(B2) \ $ C_{d}^{\prime} ( \Q _{p}) \neq \emptyset $ and $
C_{d}^{\prime} ( \Q _{q}) \neq \emptyset . $
\par  \vskip 0.1 cm
(B3) \ For each prime $ l\mid D , l \nmid d ,\quad C_{d}^{\prime} (
\Q _{l}) \neq \emptyset \Longleftrightarrow \ \left( 1 -
\left(\frac{d}{l}\right) \right) \left( 1 - \left(\frac{p q
d}{l}\right) \right) = 0. $
\par  \vskip 0.2 cm
(B4) \ For each prime $ l \mid D , l \mid d, \ C_{d}^{\prime} ( \Q
_{l}) \neq \emptyset \Longleftrightarrow \left( 1 - \left(
\frac{pdDl^{-2}}{l}\right) \right) \left( 1 -
\left(\frac{qdDl^{-2}}{l}\right) \right)=0. $
\par  \vskip 0.2 cm

Now let $ G = (V, E) $ be a directed graph. Recall that a partition
$ (V_{1}, V_{2}) $ of $ V $ is called even if for any vertex, $ P
\in V_{2} (V_{1}), \sharp \{ P \rightarrow V_{1}(V_{2}) \} $ is
even. In this case, we shall write $ (V_{1}, V_{2}) \mapsto_{e} V. $
The partition $ (V_{1}, V_{2}) $ is called quasi-even if for any
vertex, $ P \in V_{1} (V_{2}), $ \\
$ \sharp \{ P \rightarrow V_{2}(V_{1}) \}  \equiv \  \left \{
   \begin{array}{l}
  0 ( \text{mod} 2) \qquad \quad \text{if} \left(\frac{ 2 }{P}\right) = 1, \\
  1 ( \text{mod} 2) \qquad \quad \text{if} \left(\frac{ 2 }{P}\right) = - 1.
  \end{array}
  \right.   $ \\
In this case, we shall write $ (V_{1}, V_{2})\mapsto_{qe} V $ (see
[F2] and [FJ] for these definitions and related facts). Throughout
this paper, for convenience, we write empty product as $ 1. $
\par  \vskip 0.2 cm

{\bf Definition 2.5}  \ Let $ D =
 D_{1}D_{2} \cdots D_{s}D_{s+1} \cdots D_{n} $
with $ \left(\frac{p q}{D_{i}}\right) = 1 \ (i \leq s) \ $ and $
\left(\frac{p q}{D_{j}}\right) = -1 \ (s < j \leq n ) $ for some
non-negative integer $ s \leq n. $ A directed graph G(+D) is defined
as follows:
\par  \vskip 0.1 cm
Case 1. \ If $ m, p $ and $ D $ satisfy one of the following
conditions: \\ (1) \ $ m = 1; $ \ (2) \ $  m = 2, (  p + 2 ) D \not
\equiv 5 ( \text{mod} 8); $ \ (3) \ $ m \geq 3,  p  D \equiv 1 (
\text{mod} 4), $ then define the directed graph $ G(+D) = G_{1}(+D)
$ by defining the vertex $ V(G(+D)) $ to be $ V(G_{1}(+D)) = \{ -1,
p, q, D_{1}, D_{2}, \cdots , D_{n} \} $  \ and the edges \\
$ E(G(+D)) $ as \ $ E(G_{1}(+D)) = \{ \overrightarrow{D_{i}D_{j}}: \
\left(\frac{D_{j}}{D_{i}}\right) = -1, 1 \leq i \leq s, 1 \leq j
\leq n  \} \ \bigcup \ \{ \overrightarrow{D_{j}D_{i}}: \
\left(\frac{D_{i}}{D_{j}}\right) = -1, 1 \leq i \leq s, s < j \leq n
\} \ \bigcup \ \{ \overrightarrow{l D_{i}}: \
\left(\frac{D_{i}}{l}\right) = -1, 1 \leq i \leq s, l = p ,q \} \
\bigcup \ \{ \overrightarrow{-1 D_{i}}: \ \left(\frac{-1}{D_{i}}
\right) = -1, \ 1 \leq i \leq s \} \ \bigcup \ \{
\overrightarrow{D_{i} p}: \ \left(\frac{ p }{D_{i}} \right) = -1, \
1 \leq i \leq s \}. $
\par  \vskip 0.1 cm
Case 2. \ If $ m, p $ and $ D $ satisfy one of the following
conditions: \\
(1) \ $  m = 2, (  p + 2 ) D \equiv 5 ( \text{mod} 8); $ \ (2) \ $ m
\geq 3,  p  D \equiv 3 ( \text{mod} 4), $ then define the directed
graph $ G(+D) = G_{2}(+D) $ by defining the vertex $ V(G(+D)) $ to
be $ V(G_{2}(+D)) = \{ p, q, D_{1}, D_{2}, \cdots , D_{n} \} $ \ and
the edges  $ E(G(+D)) $ as \\
$ E(G_{5}^{2}(+D)) = \{ \overrightarrow{D_{i}D_{j}}: \
\left(\frac{D_{j}}{D_{i}} \right) = -1, \ 1 \leq i \leq s, 1 \leq j
\leq n  \} \ \bigcup \ \{ \overrightarrow{D_{j}D_{i}}: \
\left(\frac{D_{i}}{D_{j}} \right) = -1, \ 1 \leq i \leq s, \ s < j
\leq n \} \ \bigcup \ \{ \overrightarrow{l D_{i}}: \
\left(\frac{D_{i}}{l} \right) = -1, \ 1 \leq i \leq s, \ l = p, q \}
\ \bigcup \ \{ \overrightarrow{D_{i} p}: \ \left(\frac{ p
}{D_{i}}\right) = -1, \ 1 \leq i \leq s \}. $ \\
Here we define $ \left(\frac{ 2 }{ -1 }\right) = 1, $
if $ m, p $ and $ D $ satisfy one of the following conditions: \\
(1) \ $ m = 1, \ p D \equiv 7 ( \text{mod} 8) $ and $ D \equiv 1 (
\text{mod} 4); $ \\
(2) \ $ m = 1, \ p D \equiv 1 ( \text{mod} 8) $ and $ D \equiv 3 (
\text{mod} 4); $ \ (3) \ $ m \geq 4, p D \equiv 1 ( \text{mod} 8). $
\\
And we define $ \left(\frac{ 2 }{ -1 }\right) = -1, $
if $ m, p $ and $ D $ satisfy one of the following conditions: \\
(1) \ $ m = 1, \ p D \equiv 5 ( \text{mod} 8) $ and $ D \equiv 1 (
\text{mod} 4); $ \\
(2) \ $ m = 1, \ p D \equiv 7 ( \text{mod} 8) $ and $ D \equiv 3 (
\text{mod} 4); $ \ (3) \ $ m \geq 4, p D \equiv 5 ( \text{mod} 8). $
\par \vskip 0.2 cm

{\bf Lemma 2.6.} \ For every even partition $ ( V_{1}, V_{2} ) $ of
$ G(+D) $ \ such that $ V_{1} $ \ contains no $ -1, p, q $ or $
D_{k} \ (s < k \leq n), $ \ we have $ d \in S^{(\varphi )} (E / \Q )
$ \ where $ d = \prod_{P_{0}\in V_{1}}P_{0}. $ \ Conversely, suppose
d is odd and $ d \in S^{(\varphi )} (E / \Q ), $ we may write $ d =
P_{1} P_{2} \cdots P_{t} $ with $ 1 \leq t \leq s $ for distinct $
P_{j} \in V(G(+D)) \ (1 \leq j \leq t), $ then $ ( V_{1}, V_{2} ) $
is even, where $ V_{1} = \{P_{1}, P_{2}, \cdots , P_{t} \}. $
\par  \vskip 0.1 cm

{\bf Proof.} \ Suppose $ ( V_{1}, V_{2} ) $ is a nontrivial even
partition of $ G(+D) $ such that $ -1, p, q, D_{k} \notin V_{1} \ (s
< k \leq n). $ Let $ V_{1} = \{ D_{1}, D_{1}, \cdots ,D_{t} \} $ for
some $ 1 \leq t \leq s. $  Consider $ d = D_{1}D_{2} \cdots D_{t}. $
For any $ 1 \leq i \leq t, $  we have $ \left(\frac{p d D
D_{i}^{-2}}{D_{i}}\right) = ( - 1 )^{\sharp
\{{\overrightarrow{D_{i}P}: P \in V_{2}} \}} = 1 $ since $ ( V_{1},
V_{2} ) $ is even. Therefore, $ C_{d} ( \Q _{D_{i}}) \neq \emptyset
$ by Proposition 2.1(C)(3). Also, for $ P \in V_{2}, P \neq -1,
\left(\frac{d}{P}\right) = ( - 1 )^{\sharp \{{\overrightarrow{P
D_{i}}: \ D_{i} \in V_{1}} \}} = 1 $ since $( V_{1}, V_{2} )$ is
even. Therefore, $ C_{d} ( \Q _{P}) \neq \emptyset $ by Proposition
2.1(C)(2). We claim that $ C_{d} ( \Q _{2}) \neq \emptyset $ since $
( V_{1}, V_{2} ) $ is even. For an example in case 1, $ m = 1: $
because $ \sharp \{ \overrightarrow{-1 D_{i}}: \ D_{i}\in V_{1} \} $
is even, $ d \equiv 1 \text{mod} 4. $ Therefore, $ C_{d} ( \Q _{2})
\neq \emptyset $ by Proposition 2.1(C)(1). The remaining cases can
be done similarly. And by Proposition 2.1(A), \ we have d in $
S^{(\varphi )} (E / \Q ). $ \\
Conversely, suppose $ d = P_{1}P_{2} \cdots P_{t} \in S^{(\varphi )}
(E / \Q ) $ and $ d $ is odd. By Proposition 2.1(C), $ P_{i} \in \{
D_{1}, D_{2}, \cdots , D_{s} \} $ and $ \left(\frac{ p d D
P_{i}^{-2}}{P_{i}} \right) = 1 $ for each $ 1 \leq i \leq t. $ Let $
V_{1} = \{P_{1}, P_{2}, \cdots , P_{t} \}. $ Therefore, $ 1 =
\left(\frac{ p d D P_{i}^{-2}}{P_{i}} \right) = (-1)^{\sharp
\{{\overrightarrow{P_{i}P}: \ P \in V_{2}}\}} $ for $ 1 \leq i \leq
t. $ So we get $ \sharp \{{\overrightarrow{P_{i}P}: P \in V_{2}} \}
$ is even. For prime $ P \mid p q D d^{-1}, $ we have $  P \in V_{2}
$ and $ \left(\frac{ d }{P} \right) = 1. $ Therefore, $ 1 =
\left(\frac{ d }{ P } \right) = (-1)^{\sharp \{{\overrightarrow{P
P_{i}}: \ P_{i} \in V_{1}}\}}, $  which shows that $ \sharp
\{{\overrightarrow{ P P_{i}}: \ P_{i} \in V_{1}} \} $ is even. If $
-1 \in V_{2} $ in case 1, then $ d \equiv 1 \text{mod} 4 $ for $
C_{d} ( \Q _{2}) \neq \emptyset . $ Hence $ \sharp  \{
\overrightarrow{-1 P_{i}}: \ 1 \leq i \leq t \} $ is even. To sum
up, $ ( V_{1}, V_{2}) $ is even. The proof of lemma 2.6 is
completed. \quad $ \Box $
\par  \vskip 0.2 cm

{\bf Lemma 2.7.} \ For every quasi-even partition $( V_{1}, V_{2} )
$ of $ G(+D) $ \ such that $ V_{1} $ \ contains no $ -1, p, q $ or $
D_{k} \ (s < k \leq n), $ \ we have $ 2 d \in S^{(\varphi )} (E / \Q
), $ \ where $ d = \prod_{P_{0} \in V_{1}}P_{0}. $ \ Conversely, If
$ d $ is even and $ d \in S^{(\varphi )} (E / \Q ), $ we may write $
d = 2 P_{1}P_{2} \cdots P_{t} $ with $ 1 \leq t \leq s $ for
distinct $ P_{j} \in V(G(+D)) \ (1 \leq j \leq t), $ then $ ( V_{1},
V_{2} ) $ is quasi-even, where $ V_{1} = \{P_{1}, P_{2}, \cdots ,
P_{t} \}. $
\par  \vskip 0.1 cm

{\bf Proof.} \ Suppose $ ( V_{1}, V_{2} ) $ is a nontrivial
quasi-even partition of $ G(+D) $ such that $ -1, p, q, D_{k} \notin
V_{1} \ (s < k \leq n). $ Let $ V_{1} = \{ D_{1}, D_{1}, \cdots ,
D_{t} \} $ for some $ 1 \leq t \leq s. $ Consider $ 2 d = 2
D_{1}D_{2} \cdots D_{t}. $ For any $ 1 \leq i \leq t, $ we have $
\left(\frac{2 p d D D_{i}^{-2}}{D_{i}} \right) = \left(\frac{ 2
}{D_{i}}\right)( - 1 )^{\sharp \{{\overrightarrow{D_{i}P}: \ P \in
V_{2}} \}} = 1 $ since $( V_{1}, V_{2} )$ is quasi-even. Therefore,
$ C_{2d} ( \Q _{D_{i}}) \neq \emptyset $ by Proposition 2.1(B)(3).
Also, for $ P \in V_{2} $ and $ P \neq -1, \ \left(\frac{ 2 d }{ P
}\right) = \left(\frac{ 2 }{P}\right) ( -1 )^{\sharp
\{{\overrightarrow{P D_{i}}: D_{i} \in V_{1}} \}} = 1 $ since $ (
V_{1}, V_{2} ) $ is quasi-even. Therefore, $ C_{2d} ( \Q _{P}) \neq
\emptyset $ by Proposition 2.1(B)(2). We assert that $ C_{2d} ( \Q
_{2}) \neq \emptyset . $ To see this, we only need to prove case 1
with $ m = 1, D \equiv 1 ( \text{mod} 4) $ \ and $ p D \equiv 7 (
\text{mod} 8), $ the other cases can be similarly done. Firstly,
since $ \left(\frac{ 2 }{-1}\right) = 1, $ we have $ \sharp \{
\overrightarrow{-1 D_{i}}: \ 1 \leq i \leq t \} $ is even. So $ d
\equiv 1 \text{mod} 4 $ and $ 2 D ( 2 d )^{-1} \equiv 1 \text{mod}
4. $ Next, since $ p D \equiv 7 ( \text{mod} 8), $ we have $ d (1 -
2 D ( 2 d )^{-1})^{2} - 2 p D \equiv 2( \text{mod} 16), $ i.e., $ d
- 2 D ( p + 1 ) + \frac{2  D^{2} }{ 2 d } \equiv 2( \text{mod} 16),
$ which shows that  $ C_{2d} ( \Q _{2}) \neq \emptyset $ by
Proposition 2.1(B)(1). Furthermore by Proposition 2.1(A), we get $
2 d \in S^{(\varphi )} (E / \Q ). $ \\
Conversely, suppose $ d = 2 P_{1}P_{2}\cdots P_{t} \in S^{(\varphi
)} (E / \Q ). $ By Proposition 2.1(B), $ P_{i} \in \{ D_{1}, D_{2}
,\cdots , D_{s} \} $ and $ \left(\frac{ p d D
P_{i}^{-2}}{P_{i}}\right) = 1 $ for each $ 1 \leq i \leq t. $  Let $
V_{1} = \{P_{1}, P_{2}, \cdots ,P_{t}\}. $ Therefore, $ 1 =
\left(\frac{ p d D P_{i}^{-2}}{P_{i}}\right) =  \left(\frac{ 2
}{P_{i}}\right) (-1)^{\sharp \{{\overrightarrow{P_{i}P}: \ P \in
V_{2}}\}} $ for $ 1 \leq i \leq t . $ So $ \sharp \{ P_{i}
\rightarrow V_{2} \} = 0 ( \text{mod} 2 ), $ if  $ \left(\frac{ 2
}{P_{i}}\right) = 1 $ \ or $ 1 ( \text{mod} 2), $ if $ \left(\frac{
2 }{P_{i}}\right) = -1. $ For prime $ P \mid 2 p q D d^{-1}, $ we
have $ P \in V_{2} $ and $ \left(\frac{ d }{P}\right) = 1. $
Therefore, $ 1 = \left(\frac{ d }{ P }\right) =  \left(\frac{ 2
}{P_{i}}\right)(-1)^{\sharp \{{\overrightarrow{P P_{i}}: \ P_{i} \in
V_{1}}\}}, $ which shows that $ \sharp  \{ P \rightarrow V_{1} \} =
0 ( \text{mod} 2), $ if $ \left(\frac{ 2 }{ P }\right) = 1 $ \ or $
1 ( \text{mod} 2), $ if $ \left(\frac{ 2 }{ P }\right) = -1. $ If $
-1 \in V_{2} $ in case 1, e.g., $ m = 1, \ D \equiv 1 ( \text{mod}
4) $ and $ p D \equiv 7 ( \text{mod} 8): $ for $ C_{2d} ( \Q _{2})
\neq \emptyset, $ by Proposition 2.1(A) we have $ d \equiv 1 (
\text{mod} 4). $ Hence $ \sharp \{ \overrightarrow{-1 P_{i}}: \ 1
\leq i \leq t \} $ is even (Here notice that $ \left(\frac{ 2
}{-1}\right) = 1 $). The remaining cases can be done similarly. To
sum up, $( V_{1}, V_{2})$ is quasi-even. The proof of lemma 2.7 is
completed. \quad $ \Box $
\par  \vskip 0.2 cm

{\bf Proof of Theorem 1.1. } By Proposition 2.1, \ $ S^{(\varphi )}
(E / \Q ) \subset \{ 2, D_{1}, D_{2}, \cdots , D_{n}\}. $
Furthermore, by lemma 2.6 and lemma 2.7, it is easy to obtain all
the corresponding results for different $ m, p, D. $  The proof is
completed. \quad $ \Box $
\par  \vskip 0.2 cm

{\bf Definition 2.8.} \ Let $ D = D_{1}D_{2} \cdots D_{s}D_{s+1}
\cdots D_{n} $ with $ \left(\frac{p q}{D_{i}}\right) = 1 \ (i \leq
s) \ $ and $ \left(\frac{p q}{D_{j}}\right) = -1 \ (s < j \leq n) $
for some non-negative integer $ s \leq n. $ A graph directed g(+D)
is defined as follows :
\par  \vskip 0.1 cm
Case 1. \ If $ m, p $ and $ D $ satisfy one of the following
conditions: \\
(1) \ $ m = 1, p \equiv 1 ( \text{mod} 4) $ and $ p - D \equiv 0, 6
( \text{mod} 8); $ \\
(2) \ $ m = 1, p \equiv 3 (
\text{mod} 4) $ and $ p - D \equiv 2, 4 ( \text{mod} 8); $ \\
(3) \ $ m = 2, p D \equiv 1 ( \text{mod} 4); $ (4) \ $ m = 2,  D
\equiv 3 ( \text{mod} 4) $ and $ p D \equiv 3 ( \text{mod} 8); $
\\
(5) \ $ m = 2, D \equiv 1 ( \text{mod} 4) $ and $ p D \equiv 7 (
\text{mod} 8); $ \ (6) \ $ m = 3, p D \equiv 1 ( \text{mod} 4 ), $
\\
then define the directed graph $ g(+D) = g_{1}(+D) $ by defining
the vertex $ V(g(+D)) $ to be $ V(g_{1}(+D)) = \{ -1, p, D_{1},
D_{2}, \cdots ,D_{n} \} $ and the edges $ E(g(+D)) $ as \\
$ E(g_{1}(+D)) = \{
\overrightarrow{D_{i}D_{j}}:\left(\frac{D_{j}}{D_{i}} \right) = -1,
1 \leq i \leq s, 1 \leq j \leq n  \} \ \bigcup \ \{ \overrightarrow{
D_{i} -1 }: \ \left(\frac{-1}{D_{i}}\right) = -1, 1 \leq i \leq s \}
\ \bigcup \ \{ \overrightarrow{D_{i} p}: \ \left(\frac{ p
}{D_{i}}\right) = -1, 1 \leq i \leq s \}. $
\par  \vskip 0.1 cm
Case 2. \ If $ m, p $ and $ D $ satisfy one of the following
conditions: \\
(1) \ $ m = 1,  p \equiv 1 ( \text{mod} 4) $ and $ p - D \equiv 2,
4( \text{mod} 8); $ \ (2) \ $ m \geq 4, p D \equiv 5 ( \text{mod}
8), $ \\
then define the directed graph $ g(+D) = g_{2}(+D) $ by
defining the vertex $ V(g(+D)) $ to be $ V(g_{2}(+D))
= \{ -1, -2, p, D_{1}, D_{2}, \cdots ,D_{n} \} $
and the edges $ E(g(+D)) $ as \\
$ E(g_{2}(+D)) = \{ \overrightarrow{D_{i}D_{j}}: \
\left(\frac{D_{j}}{D_{i}}\right) = -1, 1 \leq i \leq s, 1 \leq j
\leq n  \} \ \bigcup \ \{ \overrightarrow{ D_{i} -1 }: \
\left(\frac{-1}{D_{i}}\right) = -1, 1 \leq i \leq s \} \ \bigcup \
\{ \overrightarrow{D_{i} p}:\left(\frac{ p }{D_{i}}\right) = -1, 1
\leq i \leq s \} \ \bigcup \ \{ \overrightarrow{ -2 D_{k} }: \
\left(\frac{-2}{D_{k}}\right) = -1, 1 \leq k \leq n \} \ \bigcup \
\{ \overrightarrow{-2 p }:\left(\frac{-2}{p} \right) = -1 \ \bigcup
\ \{ \overrightarrow{-2 -1 } \}. $
\par  \vskip 0.1 cm
Case 3. \ If $ m, p $ and $ D $ satisfy one of the following
conditions: \\
(1) \ $ m = 1, p \equiv 3 ( \text{mod} 4) $ and $ p - D \equiv 0, 6(
\text{mod} 8); $ \ (2) \ $ m \geq 4, p D \equiv 1 ( \text{mod} 8), $
\\
then define the directed graph $ g(+D) = g_{3}(+D) $ by defining
the vertex $ V(g(+D)) $ to be $ V(g_{3}(+D)) = \{ -1 ,
p, 2, D_{1}, D_{2}, \cdots ,D_{n} \} $ \ and the edges $ E(g(+D)) $ as \\
$ E(g_{3}(+D)) = \{
\overrightarrow{D_{i}D_{j}}:\left(\frac{D_{j}}{D_{i}}\right) = -1, 1
\leq i \leq s, 1 \leq j \leq n  \}  \ \bigcup \ \{ \overrightarrow{
D_{i} -1}: \ \left(\frac{-1}{D_{i}}\right) = -1, 1 \leq i \leq s \}
\ \bigcup \ \{ \overrightarrow{D_{i} p}: \ \left(\frac{ p
}{D_{i}}\right) = -1, 1 \leq i \leq s \} \ \bigcup \ \{
\overrightarrow{ 2 D_{k} }:\left(\frac{2}{D_{k}} \right) = -1, 1
\leq k \leq n \} \ \bigcup \ \{ \overrightarrow{2 p }: \
\left(\frac{2}{p}\right) = -1 \}. $
\par  \vskip 0.1 cm
Case 4. \ If $ m, p $ and $ D $ satisfy one of the following
conditions: \\
(1) \ $ m = 2, D \equiv 1 ( \text{mod} 4) $ and $ p D \equiv 3 (
\text{mod} 8); $ \\
(2) \ $ m = 2, D \equiv 3 ( \text{mod} 4) $ and $  p D \equiv 7 (
\text{mod} 8); $ \\
(3) \ $ m = 3, p D \equiv 3 ( \text{mod} 8); $ \ (4) \ $ m = 4, p D
\equiv 3 ( \text{mod} 4); $ \ (5) \ $ m \geq 5, p D \equiv 3 (
\text{mod} 8), $ \\
then define the directed graph $ g(+D) =
g_{4}(+D) $ by defining the vertex $ V(g(+D)) $ to be $ V(g_{4}(+D))
= \{ -1, p, D_{1}, D_{2}, \cdots , D_{n} \} $
and the edges $ E(g(+D)) $ as \\
$ E(g_{4}(+D)) = \{
\overrightarrow{D_{i}D_{j}}:\left(\frac{D_{j}}{D_{i}} \right) = -1 ,
1 \leq i \leq s, 1 \leq j \leq n  \} \ \bigcup \ \{ \overrightarrow{
D_{i} -1 }: \ \left(\frac{-1}{D_{i}}\right) = -1 , 1 \leq i \leq s
\} \ \bigcup \ \{ \overrightarrow{D_{i} p}: \ \left(\frac{ p
}{D_{i}}\right) = -1 , 1 \leq i \leq s \} \ \bigcup \ \{
\overrightarrow{-1 D_{k}}: \ \left(\frac{-1}{D_{k}}\right) = -1, 1
\leq k \leq n \} \ \bigcup  \ \{ \overrightarrow{-1 p}: \left(\frac{
-1 }{ p }\right) = -1 \}. $
\par  \vskip 0.1 cm

Case 5. \ If $ m, p $ and $ D $ satisfy one of the following
conditions: \\
(1) \ $ m = 3, p D \equiv 7 ( \text{mod} 8 ); $ \ (2) \ $ m \geq 5,
p D \equiv 7 ( \text{mod} 8), $ \\
then define the directed graph $ g(+D) = g_{5}(+D) $ by defining the
vertex $ V(g(+D)) $ to be $ V(g_{5}(+D)) = \{ -1, p, -2, 2,
D_{1}, D_{2}, \cdots , D_{n} \} $ and the edges $ E(g(+D)) $ as \\
$ E(g_{5}(+D)) = \{ \overrightarrow{D_{i}D_{j}}: \
\left(\frac{D_{j}}{D_{i}}\right) = -1, 1 \leq i \leq s, 1 \leq j
\leq n  \} \ \bigcup \ \{ \overrightarrow{ D_{i} -1}: \
\left(\frac{-1}{D_{i}}\right) = -1, 1 \leq i \leq s \} \ \bigcup \
\{ \overrightarrow{D_{i} p}: \ \left(\frac{ p }{D_{i}}\right) = -1,
1 \leq i \leq s \} \ \bigcup \ \{ \overrightarrow{ -2 D_{k} }: \
\left(\frac{-2}{D_{k}}\right) = -1, 1 \leq k \leq n \} \ \bigcup \
\{ \overrightarrow{-2 p}: \ \left(\frac{-2}{p}\right) = -1 \} \
\bigcup \ \{ \overrightarrow{ 2 D_{k} }: \
\left(\frac{2}{D_{k}}\right) = -1, 1 \leq k \leq n \} \ \bigcup \ \{
\overrightarrow{ 2 p }: \ \left(\frac{2}{p}\right) = -1 \} \ \bigcup
\ \{ \overrightarrow{-2 -1 } \}. $
\par  \vskip 0.2 cm

{\bf Lemma 2.9.} \ For every even partition $ ( V_{1}, V_{2} ) $ of
$ g(+D) $ such that $ V_{1} $ contains no $ \pm 2, $ we have $ d \in
S^{(\widehat{\varphi })} (E^{\prime } / \Q ), $ where $ d =
\prod_{P_{0} \in V_{1}}P_{0}. $ \ Conversely, if $ d $ is odd and $
d \in S^{(\widehat{\varphi })} (E^{\prime } / \Q ), $ we may write $
d = P_{1}P_{2} \cdots P_{t} $ for distinct $ P_{j} \in V(g(+D)) \ (1
\leq j \leq t), $ then $ ( V_{1}, V_{2} ) $ is even, where $ V_{1} =
\{P_{1}, P_{2}, \cdots , P_{t} \}. $
\par  \vskip 0.1 cm
{\bf Proof.} \ Suppose $ ( V_{1}, V_{2} ) $ is a nontrivial even
partition of $ g(+D) $ such that $ \pm 2 \notin V_{1}. $ Let $ V_{1}
= \{ P_{1}, P_{1}, \cdots , P_{t} \}, P_{i} \in \{-1, p, D_{1},
D_{2}, \cdots, D_{n} \} $ for each $ 1 \leq i \leq s. $ Consider $ d
= P_{1}P_{2} \cdots P_{t}. $ For each prime $ l \mid \gcd (D, \ d),
$ if $ l \in \{ D_{j}: \ s < j \leq n \}, \ (\left(\frac{-p d D
l^{-2}}{l}\right)-1) ( \left(\frac{-q d D l^{-2}}{l}\right)-1) = 0 $
because $\left(\frac{-p d D l^{-2}}{l}\right) \left(\frac{-q d D
l^{-2}}{l}\right) = \left(\frac{p q}{l}\right)  = -1; $ if $ l \in
\{ D_{i}: \ 1 \leq i \leq s \}, $ then $ \left(\frac{-p d D
l^{-2}}{l} \right) = ( - 1 )^{\sharp \{{\overrightarrow{l P}: \ P
\in V_{2}}\}} = 1 $ because $ (V_{1},V_{2}) $ is even. Therefore, by
Proposition2.2(B)(B4), we have $ C_{d}^{\prime} ( \Q_{l} ) \neq
\emptyset. $ Also for each prime $ l $ such that $ l \mid D $ and $
l \nmid d, $ if $ l \in \{ D_{j}: \ s < j \leq n \}, $ then $
(\left(\frac{ d }{l}\right)-1) ( \left(\frac{p q d }{l}\right)-1) =
0 $ because $ \left(\frac{ d }{l}\right) \left(\frac{p q
d}{l}\right) = \left(\frac{p q}{l}\right)  = -1; $ if $ l \in \{
D_{i}: \ 1 \leq i \leq s \}, $ then $ \left(\frac{ d }{ l } \right)
= ( - 1 )^{\sharp \{{\overrightarrow{l P}: \ P \in V_{1}}\}} = 1 $
because $ (V_{1}, V_{2}) $ is even. So by Proposition2.2(B)(B3), we
have $ C_{d}^{\prime} ( \Q_{l} ) \neq \emptyset . $ We assert that $
C_{d}^{\prime} ( \Q_{2} ) \neq \emptyset . $ To see this, we only
need to prove the case 3 with $ m = 1,  p \equiv 3 ( \text{mod} 4) $
and $ p - D \equiv 0, 6( \text{mod} 8), $ the other cases can be
similarly done. In fact, since $ 2 \in V_{2} $ and $ \sharp \{
\overrightarrow{2P}: \ P \in V_{1}\} $ is even, we have $ d \equiv
\pm 1 (\text{mod} 8). $ So by Proposition2.2(B)(B2), $
C_{d}^{\prime} ( \Q_{2} ) \neq \emptyset . $ This proves our
assertion. So by Proposition2.2(B)(B3) and (A)(2), we obtain that $
d \in
S^{(\widehat{\varphi })} (E^{\prime } / \Q ). $ \\
Conversely, suppose $ d = P_{1} P_{2} \cdots P_{t} \in
S^{(\widehat{\varphi })} (E^{\prime } / \Q ) $ with distinct $ \
P_{1}, \cdots , P_{t} \in \{-1, p, D_{1}, D_{2}, \cdots, D_{n} \}. $
Let \ $ V_{1} = \{P_{1}, P_{2}, \cdots , P_{t} \}. $ For each prime
$ l $ satisfying $ l \mid \gcd (D, \ d), $ if $ l \in \{ D_{j}: \ s
< j \leq n \}, \ \sharp \{ \overrightarrow{lP} : P \in V_{2} \} $ is
even because $ \sharp \{ \overrightarrow{lP}: \ P \in V_{2} \} = 0.
$ If $ l \in \{ D_{i}: \ 1 \leq i \leq s \}, $ by
Proposition2.2(B)(B4) and $ \left(\frac{pq}{l}\right) = 1, $ we have
$ \left(\frac{-p d D l^{-2}}{l} \right)-1 = 0, $ and so $ 1 =
\left(\frac{-p d D l^{-2}}{l}\right) = ( - 1 )^{\sharp
\{{\overrightarrow{lP}: P \in V_{2}}\}}, $ which shows that $
{\sharp \{{\overrightarrow{lP}: \ P \in V_{2}}\}} $ is even. Also,
for each prime $ l $ satisfying $ l \mid D $ and $ l \nmid d, $ if $
l \in \{ D_{j}: \ s < j \leq n \}, $ then $ \sharp \{
\overrightarrow{lP}: P \in V_{2} \} $ is even because $ \sharp \{
\overrightarrow{lP}: P \in V_{2} \} = 0; $ if $ l \in \{ D_{i}: \ 1
\leq i \leq s \}, $ by Proposition2.2(B)(B3) and $
\left(\frac{pq}{l}\right) = 1, $ we have $ \left(\frac{d}{l} \right)
-1 = 0. $ So $ 1 = \left(\frac{d}{l} \right) = ( - 1 )^{\sharp
\{{\overrightarrow{lP}: \ P \in V_{1}}\}}, $ which shows that $
{\sharp \{{\overrightarrow{lP}: P \in V_{1}}\}}  $ is even. As for
the vertex $ l = p, -1, $ by the definition of $ g(+D), $ we have $
\sharp \{ \overrightarrow{ l P}: \ P \in V_{1} \} = 0 $ or $ \sharp
\{ \overrightarrow{ l P}: P \in V_{2} \} = 0. $ Now firstly, in case
2, $ -2 \in V_{2}. $ By $ C_{d}^{\prime} ( \Q_{2} ) \neq \emptyset $
and the conditions for $ m, p, D, $ we have $ d \equiv 1, 3
(\text{mod} 8 ). $ So \ $ {\sharp \{{\overrightarrow{-2P}: \ P \in
V_{1}}\}} $ is even. Secondly, in case 3, $ 2 \in V_{2}. $ By $
C_{d}^{\prime} ( \Q_{2} ) \neq \emptyset $ and the conditions for $
m, p, D, $ we have $ d \equiv 1, 7 (\text{mod} 8 ). $ So \ $ {\sharp
\{{\overrightarrow{2P}: \ P \in V_{1}}\}} $ is even. Lastly, in case
5, $ \pm 2 \in V_{2}. $ By  $ C_{d}^{\prime} ( \Q_{2} ) \neq
\emptyset $ and the conditions for $ m, p, D, $ we have $ d \equiv 1
(\text{mod} 8 ). $ So both $ {\sharp \{{\overrightarrow{-2P}, P \in
V_{1}}\}} $ and $ {\sharp \{{\overrightarrow{2P}: P \in V_{1}}\}} $
are even. To sum up, $( V_{1}, V_{2} )$ is even. The Proof is
completed.  \quad $ \Box $
\par  \vskip 0.2 cm

{\bf Proof of Theorem 1.2.} \ By Proposition 2.2, we have $ \{ 1,
pq, -pD, -qD \} \subset S^{(\widehat{\varphi })} (E^{\prime } / \Q
). $ Then the conclusion follows easily by Lemma2.9. The proof is
completed. \quad $ \Box $
\par  \vskip 0.2 cm

{\bf Definition 2.10.} \ Let $ D = D_{1}D_{2} \cdots D_{s}D_{s+1}
\cdots D_{n} $ with $ \left(\frac{p q}{D_{i}}\right) = 1 \ (i \leq
s) $ and $ \left(\frac{p q}{D_{j}}\right) = -1 \ (s < j \leq n) $
for some non-negative integer $ s \leq n. $ A directed graph
$ G(-D) $ is defined as follows: \\
Case 1. If $ m =1 , D \equiv 1 ( \text{mod} 4), $ then define the
directed graph $ G(-D) = G_{1}(-D) $ by defining the vertex $
V(G(-D)) $ to be $ V(G_{1}(-D)) = \{ -1 ,p, q, D_{1}, D_{2}, \cdots
,D_{n} \} $ and the edges $ E(G(-D)) $ as \ $ E(G_{1}(-D))= \{
\overrightarrow{D_{i}D_{j}}: \left(\frac{D_{j}}{D_{i}} \right) = -1
, 1 \leq i \leq s, 1 \leq j \leq n  \} \ \bigcup \ \{
\overrightarrow{D_{j}D_{i}}: \ \left(\frac{D_{i}}{D_{j}}\right) =
-1, 1 \leq i \leq s, \ s < j \leq n  \} \ \bigcup \ \{
\overrightarrow{l D_{i}}: \ \left(\frac{D_{i}}{l}\right) = -1, 1
\leq i \leq s, \ l = p, q \} \ \bigcup \ \{ \overline{ D_{k} -1}: \
\left(\frac{-1}{D_{k}}\right) = -1, 1 \leq k \leq n \} \ \bigcup \
\{ \overrightarrow{D_{i} p}: \ \left(\frac{ p }{D_{i}}\right) = -1,
1 \leq i \leq s \} \ \bigcup \ \{ \overline {-1 l }: \
\left(\frac{-1}{l} \right) = -1, l = p, q \}. $ \\
Case 2. \ If $ m, p $ and $ D $ satisfy one of the following
conditions: \\
(1) \ $ m =1, D \equiv 3 ( \text{mod} 4). $ \ (2) \ $ m = 2, D
\equiv 3 ( \text{mod} 4) $ and $ ( p + 2 ) D \not \equiv 3 (
\text{mod} 8), $ then define the directed graph $ G(-D) = G_{2}(-D)
$ by defining the vertex $ V(G(-D)) $ to be $ V(G_{2}(-D)) = \{ -1
,p, q, D_{1}, D_{2}, \cdots ,D_{n} \} $ and the edges $
E(G(-D)) $ as \\
$ E(G_{2}(-D)) = \{
\overrightarrow{D_{i}D_{j}}:\left(\frac{D_{j}}{D_{i}}\right) = -1, 1
\leq i \leq s, \ 1 \leq j \leq n \}  \ \bigcup \ \{
\overrightarrow{D_{j}D_{i}}:\left(\frac{D_{i}}{D_{j}}\right) = -1, 1
\leq i \leq s, s < j \leq n \} \ \bigcup \ \{ \overrightarrow{l
D_{i}}: \ \left(\frac{D_{i}}{l}\right) = -1, 1 \leq i \leq s, \ l =
p, q \} \bigcup \ \{ \overline{ -1 D_{k} }: \
\left(\frac{-1}{D_{k}}\right) = -1, 1 \leq k \leq n \} \ \bigcup \
\{ \overrightarrow{D_{i} p}: \ \left(\frac{ p }{D_{i}}\right) = -1,
1 \leq i \leq s \} \ \bigcup \ \{ \overrightarrow{l -1 }: \
\left(\frac{-1}{l} \right) = -1, \ l = p, q \}. $ \\
Case 3. \ If $ m, p $ and $ D $ satisfy one of the following
conditions: \\
(1) \ $ m = 2, ( p + 2 ) D \equiv 3 ( \text{mod} 8); $ \ (2) \ $ m
\geq 3, p D \equiv 1 ( \text{mod} 4), $ \\
then define the directed graph $ G(-D) = G_{3}(-D) $ by defining the
vertex $ V(G(-D)) $ to be $ V(G_{3}(-D)) = \{ -1, p, q, D_{1},
D_{2}, \cdots , D_{n} \} $ and the edges $ E(G(-D)) $ as \\
$ E(G_{3}(-D))  = \{ \overrightarrow{D_{i}D_{j}}: \
\left(\frac{D_{j}}{D_{i}}\right) = -1, 1 \leq i \leq s, 1 \leq j
\leq n  \} \ \bigcup \ \{ \overrightarrow{D_{j}D_{i}}: \
\left(\frac{D_{i}}{D_{j}}\right) = -1, 1 \leq i \leq s, \ s < j \leq
n  \} \ \bigcup \ \{ \overrightarrow{l D_{i}}: \
\left(\frac{D_{i}}{l} \right) = -1, 1 \leq i \leq s, l = p ,q \} \
\bigcup \ \{ \overrightarrow{D_{i} p}:\left(\frac{ p }{D_{i}}\right)
= -1, 1 \leq i \leq s \} \ \bigcup \ \{ \overrightarrow{ D_{k} -1 }:
\ \left(\frac{-1}{D_{k}}\right) = -1, 1 \leq k \leq n \} \ \bigcup \
\{ \overrightarrow{ l -1 }: \
\left(\frac{-1}{l} \right) = -1, l = p, q \}. $ \\
Case 4. \ If $ m = 2, ( p + 2 )D \not \equiv 3 ( \text{mod} 8) $ and
$ D \equiv 1 ( \text{mod} 4 ), $ define the directed graph $ G(-D) =
G_{4}(-D) $ by defining the vertex $ V(G(-D)) $ to be $ V(G_{4}(-D))
= \{ -1, p, q, D_{1}, D_{2}, \cdots , D_{n} \} $ and
the edges $ E(G(-D)) $ as \\
$ E(G_{4}(-D)) = \{
\overrightarrow{D_{i}D_{j}}:\left(\frac{D_{j}}{D_{i}}\right) = -1, 1
\leq i \leq s, 1 \leq j \leq n  \} \ \bigcup \ \{
\overrightarrow{D_{j}D_{i}}: \ \left(\frac{D_{i}}{D_{j}}\right) =
-1, 1 \leq i \leq s, \ s < j \leq n  \} \ \bigcup \ \{
\overrightarrow{l D_{i}}: \ \left(\frac{D_{i}}{l}\right) = -1, 1
\leq i \leq s, \ l = p, q \} \ \bigcup \ \{ \overrightarrow{D_{i}
p}: \ \left(\frac{ p }{D_{i}}\right) = -1, 1 \leq i \leq s \} \
\bigcup \ \{ \overline{ D_{k} -1 }: \ \left(\frac{-1}{D_{k}}\right)
= -1, 1 \leq k \leq n \} \ \bigcup \ \{ \overrightarrow{ l  -1 }: \
\left(\frac{-1}{l}\right) = -1, l = p, q \} \
\bigcup \ \{ \overrightarrow{ -1  p } \}. $ \\
Case 5. \ If $ m \geq 3, p D \equiv 3 ( \text{mod} 4 ), $ define the
directed graph $ G(-D) = G_{5}(-D) $ by defining the vertex $
V(G(-D)) $ to be $ V(G_{5}(-D)) = \{ -1, p, q, D_{1}, D_{2}, \cdots
,D_{n} \} $ and the edges $ E(G(-D)) $ as  \\
$ E(G_{5}(-D)) = \{
\overrightarrow{D_{i}D_{j}}:\left(\frac{D_{j}}{D_{i}}\right) = -1, 1
\leq i \leq s, 1 \leq j \leq n  \} \ \bigcup \ \{
\overrightarrow{D_{j}D_{i}}:\left(\frac{D_{i}}{D_{j}}\right) = -1, 1
\leq i \leq s, \ s < j \leq n  \} \ \bigcup \ \{
\overrightarrow{D_{i} p}: \ \left(\frac{ p }{D_{i}}\right)= -1, 1
\leq i \leq s \} \ \bigcup \ \{ \overrightarrow{l D_{i}}: \
\left(\frac{D_{i}}{l}\right) = -1, 1 \leq i \leq s, l = p, q \} \
\bigcup \ \{ \overline{ D_{k} -1 }: \ \left(\frac{-1}{D_{k}} \right)
= -1, 1 \leq k \leq n \} \ \bigcup \ \{ \overrightarrow{ l -1 }: \
\left(\frac{-1}{l}\right) = -1, l = p, q \} \ \bigcup \ \{
\overrightarrow{ -1 p }: \ \left(\frac{-1}{p}\right) = -1 \}. $ \\
Here we define $ \left(\frac{ 2 }{ -1 }\right) = 1, $
if $ m, p $ and $ D $ satisfy one of the following conditions: \\
(1) \ $ m = 1, \ p D \equiv 7 ( \text{mod} 8) $ and $ D \equiv 1 (
\text{mod} 4); $ \\
(2) \ $ m = 1, \ p D \equiv 1 ( \text{mod} 8) $ and $ D \equiv 3 (
\text{mod} 4); $ \ (3) \ $ m = 3, p D \equiv 1 ( \text{mod}
8); $ \\
(4) \ $ m \geq 4, \ p D \equiv 7 ( \text{mod} 8); $  \ (5) \ $ m
\geq 5, p D \equiv 1 ( \text{mod} 8). $ \\
And we define $ \left(\frac{ 2 }{ -1 }\right) = -1, $
if $ m, p $ and $ D $ satisfy one of the following conditions: \\
(1) \ $ m = 1, \ p D \equiv 1 ( \text{mod} 8) $ and $ D \equiv 1 (
\text{mod} 4); $ \\
(2) \ $ m = 1, \ p D \equiv 3 ( \text{mod} 8) $ and $ D \equiv 3 (
\text{mod} 4); $ \ (3) \ $ m \geq 4, p D \equiv 3 ( \text{mod} 8). $
\par  \vskip 0.2 cm

{\bf Lemma 2.11.} \ For every even partition $ ( V_{1}, V_{2} ) $ of
$ G(-D) $ such that $ V_{1} $ contains no $ p, q $ or $ D_{k} \ (s <
k \leq n), $ \ we have $ d \in S^{(\varphi )} (E / \Q ), $ where $ d
= \prod_{P_{0}\in V_{1}}P_{0}. $ \ Conversely, if $ d $ is odd and $
d \in S^{(\varphi )} (E / \Q ), $ we may write $ d = \delta
P_{1}P_{2}\cdots P_{t} $ for $ \delta = \pm 1 $ and distinct $ P_{j}
\in V(G(-D)) \ (1 \leq j \leq t ), $ then $ ( V_{1}, V_{2} ) $ is even.
Here \\
$ V_{1} = \left \{
\begin{array}{l}
\{P_{1}, P_{2}, \cdots , P_{t} \} \quad \text{if} \ \delta = 1, \\
\{-1, P_{1}, P_{2}, \cdots ,P_{t} \} \quad \text{if} \ \delta = -1.
\end{array}
\right.   $
\par  \vskip 0.1 cm
{\bf Proof.} \ Similar to the proof of Lemma 2.6.
\par  \vskip 0.2 cm

{\bf Lemma 2.12.} \ For every quasi-even partition $( V_{1}, V_{2}
)$ of $ G(-D) $ such that $ V_{1} $ contains no $ p, q $ or $ D_{k}
\ (s < k \leq n), $ we have $ 2d \in S^{(\varphi )} (E / \Q ), $
where $ d = \prod_{P_{0}\in V_{1}}P_{0}. $ \ Conversely, if $ d $ is
even and $ d \in S^{(\varphi )} (E / \Q ), $ we may write $ d = 2
\delta P_{1}P_{2}\cdots P_{t} $ for $ \delta = \pm 1 $ and distinct
$ P_{j} \in V(G(-D)) \ (1 \leq j \leq t), $ then $( V_{1}, V_{2} ) $
is quasi-even. Here \\
$ V_{1} = \left \{
\begin{array}{l}
\{P_{1}, P_{2}, \cdots , P_{t} \} \quad \text{if} \ \delta = 1, \\
\{-1, P_{1}, P_{2}, \cdots ,P_{t} \} \quad \text{if} \ \delta = -1.
\end{array}
\right.   $
\par  \vskip 0.1 cm

{\bf Proof.} \ Similar to the proof of Lemma 2.7.
\par  \vskip 0.2 cm

{\bf Definition 2.13.} \ Let $ D = D_{1}D_{2} \cdots D_{s}D_{s+1}
\cdots D_{n} $ with $ \left(\frac{p q}{D_{i}}\right) = 1 \ (i \leq
s) $ and $ \left(\frac{p q}{D_{j}}\right) = -1 \ (s < j \leq n) $
for some non-negative integer $ s \leq n. $ A graph directed
$ g(-D) $ is defined as follows: \\
Case 1. If $ m, p $ and $ D $ satisfy one of the following
conditions: \\
(1) \ $ m = 1, p - D \equiv 2, 4 ( \text{mod} 8); $ \\
(2) \ $ m = 2, p D \equiv 3 ( \text{mod} 4); $ \ (3) \  $ m = 2, D
\equiv 1( \text{mod} 4); $ and $ p D \equiv 5( \text{mod} 8); $ \\
(4) \ $ m = 2, D \equiv 3( \text{mod} 4) $ and $ p D \equiv 1(
\text{mod} 8); $ \ (5) \  $ m = 3, p D \equiv 3 ( \text{mod} 4), $ \\
then define the directed graph $ g(-D) = g_{1}(-D) $ by defining the
vertex $ V(g(-D)) $ to be $ V(g_{1}(-D)) = \{ p, D_{1}, D_{2},
\cdots , D_{n} \} $ and the edges $ E(g(-D)) $ as $ E(g_{1}(-D)) \\
= \{ \overrightarrow{D_{i}D_{j}}: \ \left(\frac{D_{j}}{D_{i}}\right)
= -1, 1 \leq i \leq s, 1 \leq j \leq n  \} \ \bigcup \ \{
\overrightarrow{D_{i} p}: \ \left(\frac{ p }{D_{i}}\right) = -1,
1 \leq i \leq s \}. $ \\
Case 2. \ If $ m, p $ and $ D $ satisfy one of the following
conditions: \\
(1) \ $ m = 1, p \equiv 1 ( \text{mod} 4) $ and $ p - D \equiv 0, 6
( \text{mod} 8); $ \ (2) \ $ m \geq 4, p D \equiv 3 ( \text{mod} 8),
$ then define the directed graph $ g(-D) = g_{2}(-D) $ by defining
the vertex $ V(g(-D)) $ to be $ V(g_{2}(-D)) = \{ p, -2, D_{1},
D_{2}, \cdots ,D_{n} \} $ and the edges $ E(g(-D)) $ as $
E(g_{2}(-D)) \\
= \{\overrightarrow{D_{i}D_{j}}: \ \left(\frac{D_{j}}{D_{i}} \right)
= -1, 1 \leq i \leq s, 1 \leq j \leq n  \} \ \bigcup \ \{
\overrightarrow{-2 p }: \ \left(\frac{-2}{p}\right) = -1 \} \
\bigcup \ \{ \overrightarrow{D_{i} p}: \ \left(\frac{ p
}{D_{i}}\right) = -1, 1 \leq i \leq s \} \ \bigcup \ \{
\overrightarrow{ -2 D_{k} }: \ \left(\frac{-2}{D_{k}}\right) = -1,
1 \leq k \leq n \} $ \\
Case 3. \ If $ m, p $ and $ D $ satisfy one of the following
conditions: \\
(1) \ $ m = 1, p \equiv 3 ( \text{mod} 4 ) $ and $ p - D \equiv 0, 6
( \text{mod} 8); $  \ (2) \ $ m \geq 5,  p D \equiv 7 ( \text{mod}
8), $ then define the directed graph $ g(-D) = g_{3}(-D) $ by
defining the vertex $ V(g(-D)) $ to be $ V(g_{3}(-D)) = \{ p, 2,
D_{1}, D_{2}, \cdots , D_{n} \} $ and the edges $ E(g(-D)) $ as
 $ E(g_{3}(-D)) \\ = \{ \overrightarrow{D_{i}D_{j}}:
\left(\frac{D_{j}}{D_{i}}\right) = -1,
 1 \leq i \leq s, 1 \leq j \leq n  \} \ \bigcup \
\{ \overrightarrow{2 p }: \ \left(\frac{2}{p}\right) = -1 \} \
\bigcup \ \{ \overrightarrow{D_{i} p}:\left(\frac{ p }{D_{i}}\right)
= -1, 1 \leq i \leq s \}
 \ \bigcup \ \{ \overrightarrow{ 2 D_{k} }: \ \left(\frac{2}{D_{k}}
 \right) = -1, 1 \leq k \leq n \}. $ \\
Case 4. \ If $ m, p $ and $ D $ satisfy one of the following
conditions: \\
(1) \ $ m = 2, D \equiv 1 ( \text{mod} 4 ) $ and $ p D \equiv 1 (
\text{mod} 8); $ \ (2) \ $ m = 2, D \equiv 3 ( \text{mod} 4) $ and $
p D \equiv 5 ( \text{mod} 8); $ \ (3) \ $ m \geq 3, p D \equiv 5 (
\text{mod} 8); $ \ (4) \ $ m = 4, p D \equiv 1 ( \text{mod} 8), $
then define the directed graph $ g(-D) = g_{4}(-D) $ by defining the
vertex $ V(g(-D)) $ to be $ V(g_{4}(-D)) = \{ p, -1,
D_{1}, D_{2}, \cdots , D_{n} \} $ and the edges $ E(g(-D)) $ as \\
$ E(g_{4}(-D)) = \{ \overrightarrow{D_{i}D_{j}}: \
\left(\frac{D_{j}}{D_{i}}\right) = -1, 1 \leq i \leq s, \ 1 \leq j
\leq n  \} \ \bigcup \ \{ \overrightarrow{-1 p }: \
\left(\frac{-1}{p}\right) = -1 \} \ \bigcup \ \{
\overrightarrow{D_{i} p}: \ \left(\frac{ p }{D_{i}}\right) = -1, 1
\leq i \leq s \} \ \bigcup \ \{ \overrightarrow{ -1 D_{k} }: \
\left(\frac{-1}{D_{k}}\right) = -1, 1 \leq k \leq n \}. $ \\
Case 5. \ If $ m, p $ and $ D $ satisfy one of the following
conditions: \\
(1) \ $ m = 3, p D \equiv 1 ( \text{mod} 8); $ \ (2) \ $ m \geq 5, p
D \equiv 1 ( \text{mod} 8 ), $ then define the directed graph $
g(-D) = g_{5}(-D) $ by defining the vertex $ V(g(-D)) $ to be $
V(g_{5}(-D)) = \{ p, -1, 2, D_{1}, D_{2}, \cdots ,D_{n} \} $ and the
edges $ E(g(-D)) $ as \\
$ E(g_{5}(-D)) = \{ \overrightarrow{D_{i}D_{j}}: \
\left(\frac{D_{j}}{D_{i}}\right) = -1, l \leq i \leq s, \ 1 \leq j
\leq n \} \ \bigcup \ \{ \overrightarrow{-1 p }: \
\left(\frac{-1}{p} \right) = -1 \} \ \bigcup \ \{
\overrightarrow{D_{i} p}: \ \left(\frac{ p }{D_{i}} \right) = -1,
 1 \leq i \leq s \} \ \bigcup \ \{ \overrightarrow{ -1 D_{k} }:
\left(\frac{-1}{D_{k}}\right) = -1, \ 1 \leq k \leq n \} \ \bigcup \
\{ \overrightarrow{ 2 D_{k} }: \ \left(\frac{2}{D_{k}} \right) = -1,
1 \leq k \leq n \} \ \bigcup \ \{ \overrightarrow{2 p }: \
\left(\frac{2}{p}\right) = -1\}. $
\par  \vskip  0.2 cm

{\bf Lemma 2.14.} \ For every even partition $ ( V_{1}, V_{2} ) $ of
$ g(-D) $ such that $ V_{1} $ contains no $ -1, \pm 2, $ we have $ d
\in  S^{(\widehat{\varphi })} (E^{\prime } / \Q ), $ where $ d =
\prod_{P_{0}\in V_{1}}P_{0}. $ \ Conversely, if $ d $ is odd and $ d
\in  S^{(\widehat{\varphi })} (E^{\prime } / \Q ), $ we may write $
d = P_{1}P_{2}\cdots P_{t} $ for distinct $ P_{j} \in V(g(-D)) \ (1
\leq j \leq t), $ then $ ( V_{1}, V_{2} ) $ is even, where $ V_{1} =
\{P_{1}, P_{2}, \cdots , P_{t} \}. $
\par  \vskip 0.1 cm

{\bf Proof.} \ Similar to the proof of Lemma 2.9.
\par  \vskip  0.2 cm

{\bf Proofs of Theorem 1.3 and 1.4.} By using Proposition 2.3, 2.4
and Lemma 2.11, 2.12, 2.14, the proofs are similar to that of
Theorem 1.1 and 1.2.  \quad $ \Box $

\par  \vskip  0.4 cm

\hspace{-0.6cm}{\bf Appendix}

\par \vskip 0.2 cm

In this appendix, by descent method, we obtain the following results
about Selmer group of the elliptic curve (1.1) for all integers $ m
\geq 2, $ which generalize the ones in [LQ] for the case $ m = 1. $
The method is the same as in [LQ] (see also [QZ] and [DW]), so we
omit the details.
\par \vskip 0.2 cm

{\bf Theorem A.1. } \ Let $ E = E_{+} $ be the elliptic curve in
(1.1) with $ \varepsilon = + 1 $ and $ l $ be an odd prime number.
For each $ i \in \{ 1, \cdots , n \}, $ denote \\
$ \Pi _{i}^{+}(D) = \delta_{i} + \left( 1 -
\left(\frac{q\widehat{D_{i}}}{D_{i}} \right) \right) \left( 1 -
\left(\frac{p\widehat{D_{i}}}{D_{i}} \right) \right) + \sum _{ l
\mid pq \widehat{D_{i}} } \left( 1 - \left(\frac{D_{i}}{l} \right)
\right), $ \\
where $ \delta_{i} = 0 $ if $ D_{i}, m, p $ and $ D $ satisfy one of
the following conditions: \\
(1) \ $ D_{i} \equiv 1 (\text{mod} \ 4), $ \ (2) \ $ m = 2, p - D
\equiv 2 (\text{mod} \ 8), $ \ (3) \ $ m \geq 3, p + D \equiv 0
(\text{mod} \ 4); $ otherwise, $ \delta_{i} = 1. $ And denote
$$ \Pi_{n+1}^{+}(D) = \delta_{n+1} + \sum _{l \mid pqD} \left(1 -
\left(\frac{2}{l}\right) \right), $$ where \ $ \delta_{n+1} = 0, $
if $ m, p $ and $ D $ satisfy one of the following conditions: \\
(1) \ $ m = 3, p D \equiv - 1 (\text{mod} \ 8), $ \ (2) \ $ m = 4, p
D \equiv 1 (\text{mod} \ 8), $ \ (3) \ $ m \geq 5; $ otherwise, $
\delta_{n+1} = 1. $ Here \ $ \left(- \right) $ is the ( Legendre )
quadratic residue symbol. And define a function $ \rho ^{+}( D ) $
by $$ \rho ^{+}( D ) = \sum _{i = 1}^{n+1} \left[ \frac{1}{1 + \Pi
_{i}^{+}(D)}\right], $$ where [x] is the greatest integer $ \leq x.
$ Then there exists a subset $ T \subset \{2, D_{1}, \cdots , D_{n}
\} $ with cardinal $ \sharp T = \rho ^{+}( D ) $ such that $
S^{(\varphi )} (E / \Q ) \supset \ < T \text{mod} (\Q^{\star^{2}})
> \ \cong \left( \Z / 2 \Z \right)^{\rho ^{+}( D )}. $ In particular,
$ \text{dim}_{2} S^{(\varphi )} (E / \Q ) \geq \rho ^{+}( D ). $
\par \vskip 0.2 cm

{\bf Theorem A.2.} \ Let $ E^{\prime} = E_{+}^{\prime} $ be the
elliptic curve in (1.2) with $ \varepsilon = + 1. $  For each $ i
\in Z(n) = \{ 1, \cdots , n \}, $ denote \\
$ \Pi _{i}^{+}(D^{\prime}) = \left( 1 - \left(\frac{-q
\widehat{D_{i}}}{D_{i}} \right) \right) \left( 1 -
\left(\frac{-p\widehat{D_{i}}}{D_{i}} \right) \right) + \sum _{j=1,
\ j \neq i} ^{n} \left( 1 - \left(\frac{D_{i}}{D_{j}} \right)
\right) \left( 1 - \left(\frac{pqD_{i}}{D_{j}} \right) \right) $ \
and \ $ \Pi _{n+1}^{+}(D^{\prime}) = \delta^{\prime }_{n+1} + \sum
_{i=1} ^{n} \left(1 - \left(\frac{-1}{D_{i}} \right) \right) \left(
1 - \left(\frac{-pq}{D_{i}} \right) \right), $ where $
\delta^{\prime }_{n+1} = 0, $ if \ $ m, p $ and $ D $ satisfy one of
the following conditions: \ (1) \ $ m = 2, p - D \not\equiv 2
(\text{mod} \ 8 ), $ \\
(2) \ $ m = 3, p - D \equiv 0 (\text{mod} \ 4), $ \ (3) \
$ m \geq 4, p - D  \equiv 0 (\text{mod} \ 8); $ \\
otherwise, $ \delta^{\prime }_{n+1} = 1. $ Here $ \left(- \right) $
is the ( Legendre ) quadratic residue symbol. \\
Take a subset $ I $ of $ Z(n) $ as follows: \\
if \ $ m = 2, $ set $ I = \{ i \in  Z(n): \ D_{i} \equiv 1
(\text{mod} 4) \} \ \bigcup \ \{ i \in  Z(n): \ D_{i} + pD \equiv 0
(\text{mod} 4) \} \ \bigcup \ \{ i \in Z(n):  D_{i} \equiv 3
(\text{mod} 4) \ \text{and} \ p - D \equiv 6 (\text{mod} 8) \}; $ \\
if $ m = 3, $ set $ I = \{ i \in  Z(n): \ D_{i} \equiv 1 (\text{mod}
8) \} \ \bigcup \ \{i \in  Z(n): \ D_{i} + pD \equiv 0 (\text{mod}
8) \} \ \bigcup \ \{ i \in  Z(n): \ D_{i} \equiv 3 (\text{mod} 4) \
\text{and} \ pD + D_{i} \equiv 4 (\text{mod} 8) \} \ \bigcup \ \{ i
\in Z(n): \ D_{i} \equiv 5 (\text{mod} 8) \
\text{and} \ pD - D_{i} \equiv 0 (\text{mod} 4) \}; $ \\
if $ m = 4, $ set $ I = \{ i \in  Z(n): \ D_{i} \equiv 1 (\text{mod}
8) \} \ \bigcup \ \{i \in Z(n): \ D_{i} + pD \equiv 0 (\text{mod} 8)
\} \ \bigcup \ \{i \in  Z(n): \ D_{i} \equiv 5 (\text{mod} 8) \
\text{and} \ pD + D_{i} \equiv 4 (\text{mod} 8) \}; $ \\
if $ m \geq 5, $ set $ I = \{ i \in Z(n): \ D_{i} \equiv 1
(\text{mod} 8) \} \ \bigcup \ \{ i \in Z(n): \ D_{i} +
pD \equiv 0 (\text{mod} 8) \}. $ \\
Define a function $ \rho ^{+}( D^{\prime} ) $ by $$ \rho ^{+}(
D^{\prime} ) = \sum _{i \in I \bigcup \{ n + 1 \}} \left[ \frac{1}{1
+ \Pi _{i}^{+}(D^{\prime })} \right], $$ where [x] is the greatest
integer $ \leq x. $ Then there exists a subset $ T \subset \{ -1,
D_{1}, \cdots , D_{n} \} $ with cardinal $ \sharp T = \rho ^{+}(
D^{\prime} ) $ such that $ S^{(\varphi )} (E / \Q ) \supset \ < T
\text{mod} (\Q^{\star^{2}}) > \ \cong \left( \Z / 2 \Z \right)^{\rho
^{+}(D^{\prime})}. $ In particular, $ \text{dim}_{2} S^{(\varphi )}
(E / \Q ) \geq \rho ^{+}( D^{\prime } ). $
\par \vskip 0.2 cm

{\bf Theorem A.3.} \ Let $ E = E_{-} $ be the elliptic curve in
(1.1) with $ \varepsilon = - 1 $ and \ $ l $ be an odd prime number.
For each $ i \in \{ 1, \cdots , n \}, $ denote \\
$ \Pi _{i}^{-}(D) = \delta_{i} + \left( 1 - \left(\frac{-q
\widehat{D_{i}}}{D_{i}} \right) \right) \left( 1 - \left(\frac{-p
\widehat{D_{i}}}{D_{i}} \right) \right) + \sum _{ l \mid pq
\widehat{D_{i}} } \left( 1 - \left(\frac{D_{i}}{l} \right) \right),
$ where $ \delta_{i} = 0, $ if $ D_{i}, m, p $ and $ D $ satisfy one
of the following conditions: \\
(1) \ $ D_{i} \equiv 1 (\text{mod} 4), $ \ (2) \ $ m = 2, p + D
\equiv 2 (\text{mod} 8), $ \ (3) \ $ m \geq 3, p - D \equiv 0
(\text{mod} 4); $ otherwise, $ \delta_{i} = 1. $ And denote
$$ \Pi_{n+1}^{-}(D) = \delta_{n+1} + \sum _{l \mid pqD} \left( 1 -
\left(\frac{2}{l} \right) \right), $$ where $ \delta_{n+1} = 0, $
if \ $ m, p $ and $ D $ satisfy one of the following conditions: \\
(1) \ $ m = 3, p D \equiv  1 (\text{mod} 8), $ (2) \ $ m = 4, p D
\equiv -1 (\text{mod} 8), $ (3) \ $ m \geq 5; $ otherwise, $
\delta_{n+1} = 1; $ and denote
$$ \Pi_{n+2}^{-}(D) = \delta_{n+2} + \sum _{l \mid pqD} \left( 1 -
\left(\frac{-1}{l}\right) \right), $$ where  $ \delta_{n + 2} = 0, $
if \ $ m, p $ and $ D $ satisfy one of the following conditions: \\
(1) \ $ p D \equiv  1 (\text{mod} 8), $ (2) \ $ m \geq 3, p D \equiv
5 (\text{mod} 8); $ otherwise, $ \delta_{n+2} = 1. $ \\
And define a function $ \rho ^{-}( D ) $ by $$ \rho ^{-}( D ) = \sum
_{i = 1} ^{n+2} \left[ \frac{1}{1 + \Pi _{i}^{-}(D)} \right], $$
where [x] is the greatest integer $ \leq x. $ Then there exists a
subset $ T \subset \{-1, 2, D_{1}, \cdots , D_{n} \} $ with cardinal
$ \sharp T = \rho ^{-}( D ) $ such that $ S^{(\varphi )} (E / \Q )
\supset \ < \{ D_{i} : \ D_{i} \in T \} \ \text{mod} \
(\Q^{\star^{2}}) > \ \cong \left( \Z / 2 \Z \right)^{\rho ^{-}( D
)}. $  In particular, $ \text{dim}_{2} S^{(\varphi )} (E / \Q ) \geq
\rho ^{-}( D ). $
\par \vskip 0.2 cm

{\bf Theorem A.4.} \ Let $ E^{\prime } = E_{-}^{\prime } $ be the
elliptic curve in (1.2) with $ \varepsilon = -1. $ For each $ i \in
Z(n) = \{ 1, \cdots , n \}, $ denote \quad \ $ \Pi
_{i}^{-}(D^{\prime }) = \\
\left(1 - \left(\frac{q \widehat{D_{i}}}{D_{i}} \right) \right)
\left(1 - \left(\frac{p \widehat{D_{i}}}{D_{i}} \right) \right) +
\sum _{j=1, \ j \neq i}^{n} \left( 1 - \left(\frac{D_{i}}{D_{j}}
\right) \right) \left( 1 - \left(\frac{pqD_{i}}{D_{j}} \right)
\right). $ Here $ \left(- \right) $ is the (Legendre) quadratic
residue symbol. Take a subset $ I $ of $ Z(n) $ as follows: \\
if \ $ m = 2, $ set $ I = \{ i \in Z(n): \ D_{i} \equiv 1
(\text{mod} 4) \} \ \bigcup \ \{i \in Z(n): \ D_{i} - pD \equiv 0
(\text{mod} 4) \} \ \bigcup \ \{ i \in Z(n): \ D_{i} \equiv 3
(\text{mod} 4) \ \text{and} \ p + D \equiv 6 (\text{mod} 8) \}; $ \\
if $ m = 3, $ set $ I = \{i \in  Z(n): \ D_{i} \equiv 1 (\text{mod}
8) \} \ \bigcup \ \{i \in Z(n): \ D_{i} - pD \equiv 0 (\text{mod} 8)
\} \ \bigcup \ \{i \in Z(n): \ D_{i} \equiv 3 (\text{mod} 4) \
\text{and} \ pD - D_{i} \equiv 4 (\text{mod} 8) \} $ \\
$ \bigcup \{ i \in Z(n):  D_{i} \equiv 5 (\text{mod} 8) \ \text{and} \
pD + D_{i} \equiv 0 (\text{mod} 4) \}; $ \\
if $ m = 4, $ set $ I = \{ i \in Z(n): \ D_{i} \equiv 1 (\text{mod}
8) \} \ \bigcup \ \{i \in  Z(n): \ D_{i} - pD \equiv 0 (\text{mod}
8) \ \bigcup \ \{i \in Z(n): \ D_{i} \equiv 5 (\text{mod} 8) \
\text{and} \ pD - D_{i} \equiv 4 (\text{mod} 8) \}; $ \\
if $ m \geq 5, $ set $ I = \{ i \in Z(n): \ D_{i} \equiv 1
(\text{mod} 8) \} \ \bigcup \ \{i \in Z(n): \ D_{i} - pD
\equiv 0 (\text{mod} 8) \}. $ \\
Define a function $ \rho ^{-}( D^{\prime} ) $ by $$ \rho ^{-}(
D^{\prime} ) = \sum _{i \in I }  \left[ \frac{1}{1 + \Pi
_{i}^{-}(D^{\prime })} \right], $$ where [x] is the greatest integer
$ \leq x. $ Then there exists a subset $ T \subset \{D_{1}, \cdots ,
D_{n} \} $ with cardinal $ \sharp T = \rho ^{-}(D^{\prime} ) $ such
that $ S^{(\varphi )} (E / \Q ) \supset \ < T \ \text{mod} \
(\Q^{\star^{2}}) > \ \cong \left( \Z / 2 \Z \right)^{\rho ^{-}(
D^{\prime} )}. $ In particular, $ \text{dim}_{2} S^{(\varphi )} (E /
\Q ) \geq \rho ^{-}( D^{\prime } ). $
\par  \vskip 0.3 cm

{ \bf Acknowledgement } \ We are grateful to Prof. Keqin Feng for
sending us his papers [F1], [F2], [FX] and other materials which are
helpful for this work.

\par  \vskip 0.3 cm
\hspace{-0.8cm} {\bf References }
\begin{description}

\item[[DW]] A. Dabrowski, M. Wieczorek, On the equation
$ y^{2} = x ( x - 2^{m}) ( x + q - 2^{m}), $ J. Number Theory, 2007,
124: 364-379.

\item[[F1]] K. Feng, Non-congruent number, odd graphs and the BSD
conjecture, Acta Arith., 1996, 80: 71-83.

\item[[F2]] K. Feng, Non-Congruent Numbers and Elliptic Curves
with Rank Zero (in Chinese), University of Science and Technology of
China Press, 2008.

\item[[FJ]] B. Faulkner, K. James, A graphical approach to
computing Selmer groups of congruent number curves, Ramanujan J.,
2007, 14: 107-129.

\item[[FX]] K. Feng, M. Xiong, On elliptic curves
$ y^{2} = x^{3} - n^{2} x $ with rank zero, J. Number Theory, 2004,
109: 1-24.

\item[[LQ]] F. Li, D. Qiu, On Several Families of Elliptic Curves
with Arbitrary Large Selmer Groups, arXiv.org: 0911.0236v1 [math.AG]
2 Nov 2009.

\item[[QZ]] D. Qiu, X. Zhang, Mordell-weil groups and selmer
groups of two types of elliptic curves, Science in China (series A),
2002, Vol.45, No.11, 1372-1380.

\item[[S]] J. Silverman, The Arithmetic of Elliptic Curves, New
York: Springer-Verlag, 1986.

\end{description}

\end{document}